\documentclass[12pt]{article}
\usepackage{amsmath,amssymb,amsfonts,amsthm}

\usepackage{graphicx}
\newenvironment{myabstract}{\par\noindent
{\bf Abstract . } \small }
{\par\vskip8pt minus3pt\rm}
\newcounter{item}[section]
\newcounter{kirshr}
\newcounter{kirsha}
\newcounter{kirshb}
\newenvironment{enumroman}{\setcounter{kirshr}{1}
\begin{list}{(\roman{kirshr})}{\usecounter{kirshr}} }{\end{list}}
\newenvironment{enumarab}{\setcounter{kirshb}{1}
\begin{list}{(\arabic{kirshb})}{\usecounter{kirshb}} }{\end{list}}
\newenvironment{athm}[1]{\vskip3mm\par\noindent
{\bf #1 }. \slshape }
{\upshape\par\vskip10pt minus3pt}
\newtheorem{theorem}{Theorem}[section]

\newtheorem{lemma}[theorem]{Lemma}
\newtheorem{corollary}[theorem]{Corollary}
\newenvironment{demo}[1]{\noindent{\bf #1.}\upshape\mdseries}
{\nopagebreak{\hfill\rule{2mm}{2mm}\nopagebreak}\par\normalfont}
\theoremstyle{definition}

\newtheorem{definition}[theorem]{Definition}

\def\Nr{{\mathfrak{Nr}}}

\def\Sg{{\mathfrak{Sg}}}
\def\Fm{{\mathfrak{Fm}}}

\def\(R)RA{{\bf (R)RA}}

\def\C{\mathbb{C}}
\def\A{{\mathfrak{A}}}
\def\B{{\mathfrak{B}}}
\def\C{{\mathfrak{C}}}

\def\PA{{\bf PS}}

\def\PA{{\bf PA}}

\def\L{{\mathfrak{L}}}

\def\Nr{{\mathfrak{Nr}}}

\def\Sg{{\mathfrak{Sg}}}
\def\Fm{{\mathfrak{Fm}}}

\def\(R)RA{{\bf (R)RA}}

\def\C{\mathbb{C}}
\def\A{{\mathfrak{A}}}
\def\B{{\mathfrak{B}}}
\def\C{{\mathfrak{C}}}

\def\M{{\mathfrak{M}}}

\def\PA{{\bf PA}}

\title{Atomic polyadic algebras of infinite dimension are completely representable} 
\author{Tarek Sayed Ahmed\\
Department of Mathematics, Faculty of Science,\\ 
Cairo University, Giza, Egypt.
  }
\begin{document}
\maketitle

\begin{myabstract}  Answering a question posed by Hodkinson, we show that for infinite ordinals $\alpha$, every atomic polyadic algebra of 
dimension $\alpha$ 
is completely representable. We infer that certain infinitary extensions of first order logic without equality enjoy a Vaught's theorem; 
atomic theories have 
atomic models. Our proof uses a neat embedding theorem together with a simple topological argument.
\footnote{ 

keywords: Algebraic Logic, polyadic algebras, complete representations}
 
\end{myabstract}

Polyadic algebras were introduced by Halmos \cite{Halmos} to provide an algebraic reflection
of the study of first order logic without equality. Later the algebras were enriched by 
diagonal elements to permit the discussion of equality. 
That the notion is indeed an adequate reflection of first order logic was 
demonstrated by Halmos' representation theorem for locally finite polyadic algebras 
(with and without equality). Daigneault and Monk  
proved a strong extension of Halmos' theorem, namely that,  
every polyadic algebra of infinite dimension (without equality) is representable \cite{DM}. 

There are several types of representations in algebraic logic. Ordinary representations are just isomorphisms from boolean algebras
with operators to a more concrete structure (having the same similarity type) 
whose elements are sets endowed with set-theoretic operations like intersection and complementation.
Complete representations, on the other hand, are representations that preserve arbitrary conjunctions whenever defined.
The notion of complete representations has turned out to be very interesting for cylindric algebras, where it is proved in \cite{HH}
that the class of completely representable algebras is not elementary. 

The correlation of atomicity to complete representations has caused a lot of confusion in the past.
It was mistakenly thought for a long time, among algebraic logicians,  that atomic representable relation and cylindric algebras
are completely representable, an error attributed to Lyndon and now referred to as Lyndon's error.

For boolean algebras, however  this is true.  The class of completely representable algebras
is simply the class of atomic ones. An analogous result holds for certain countable reducts of polyadic algebras \cite{s}.
The notion of complete representations has been linked to Martin's axiom, omitting types theorems 
nd existence of atomic models in model theory \cite{t}, \cite{ANS}, \cite{Moh4}. Such connections 
will be further elaborated upon below in a new setting.

In this paper we show that an atomic polyadic algebra of infinite dimension is also completely representable, 
in sharp contrast to the cylindric and quasipolyadic equality cases
\cite{HH}, \cite{Moh1}. This result answers a question raised by Ian Hodkinson, see p. 260 in \cite{h} and Remark 6.4, p. 283 in op.cit. 
Furthermore, we show this result has non-trivial model-theoretic consequences 
for certain infinitary logics.

The paper is organized as follows. In section 1 we prove our main result as in the title. In section 2, we provide a non-trivial logical counterpart
of our algebraic result. In the final section we comment on related results and pose some questions.

\section{ Main result}

Our notation is in conformity with \cite{HMT2} with only one deviation. We write $f\upharpoonright A$ instead of $A\upharpoonright f$, 
for the restriction of a function $f$ to 
a set $A,$ which is the more usual standard notation. 
On the other hand, following \cite{HMT2}, for given sets $A, B$ we write $A\sim B$ for the set $\{x\in A: x\notin B\}$.
Gothic letters are used for algebras, and the corresponding Roman letter will denote their universes.
Also for an algebra $\A$ and $X\subseteq A$, $\Sg^{\A}X$, or simply $\Sg X$ when $\A$ is clear from context, denotes the subalgebra of $\A$ 
generated by $X$. $Id$ denotes the identity function and when we write $Id$ we will be tacitly assuming that its domain is clear from context.
We now recall the definition of polyadic algebras as formulated in \cite{HMT2}. 
Unlike Halmos' formulation, the dimension of algebras is specified by ordinals as opposed to 
arbitrary sets.

\begin{definition} Let $\alpha$ be an ordinal. By a {\it polyadic algebra} of dimension $\alpha$,  
or a $\PA_{\alpha}$ for short,
we understand an algebra of the form
$$\A=\langle A,+,\cdot ,-,0,1,{\sf c}_{(\Gamma)},{\sf s}_{\tau} \rangle_{\Gamma\subseteq \alpha ,\tau\in {}^{\alpha}\alpha}$$
where ${\sf c}_{(\Gamma)}$ ($\Gamma\subseteq  \alpha$) and ${\sf s}_{\tau}$ ($\tau\in {}^{\alpha}\alpha)$ are unary 
operations on $A$, such that postulates  
below hold for $x,y\in A$, $\tau,\sigma\in {}^{\alpha}\alpha$ and 
$\Gamma, \Delta\subseteq \alpha$ 

\begin{enumerate}
\item $\langle A,+,\cdot, -, 0, 1\rangle$ is a boolean algebra
\item ${\sf c}_{{(\Gamma)}}0=0$

\item $x\leq {\sf c}_{{(\Gamma)}}x$

\item ${\sf c}_{(\Gamma)}(x\cdot {\sf c}_{(\Gamma)}y)={\sf c}_{(\Gamma)}x \cdot {\sf c}_{(\Gamma)}y$

\item ${\sf c}_{(\Gamma)}{\sf c}_{(\Delta)}x={\sf c}_{(\Gamma\cup \Delta)}x$

\item ${\sf s}_{\tau}$ is a boolean endomorphism

\item ${\sf s}_{Id}x=x$

\item ${\sf s}_{\sigma\circ \tau}={\sf s}_{\sigma}\circ {\sf s}_{\tau}$

\item if $\sigma\upharpoonright (\alpha\sim \Gamma)=\tau\upharpoonright (\alpha\sim \Gamma)$, then 
${\sf s}_{\sigma}{\sf c}_{(\Gamma)}x={\sf s}_{\tau}{\sf c}_{(\Gamma)}x$


\item If $\tau^{-1}\Gamma=\Delta$ and $\tau\upharpoonright \Delta $ is
one to one, then ${\sf c}_{(\Gamma)}{\sf s}_{\tau}x={\sf s}_{\tau}{\sf c}_{(\Delta )}x$.

\end{enumerate}

\end{definition}
We will sometimes add superscripts to cylindrifications and substitutions indicating the algebra they are evaluated in.
The class of representable algebras is defined via set - theoretic operations
on sets of $\alpha$-ary sequences. Let $U$ be a set. 
For $\Gamma\subseteq \alpha$ and $\tau\in {}^{\alpha}\alpha$, we set
$${\sf c}_{(\Gamma)}X=\{s\in {}^{\alpha}U: \exists t\in X,\ \forall j\notin \Gamma, t(j)=s(j) \}$$
and 
$${\sf s}_{\tau}X=\{s\in {}^{\alpha}U: s\circ \tau\in X\}.$$
For a set $X$, let $\B(X)$ be the boolean set algebra $(\wp(X), \cup, \cap, \sim).$
The class of representable polyadic algebras, or ${\bf RPA}_{\alpha}$ for short, is defined by
$$SP\{\langle \B(^{\alpha}U), {\sf c}_{(\Gamma)}, {\sf s}_{\tau} \rangle_{\Gamma\subseteq \alpha, \tau\in {}^{\alpha}\alpha}:
\ \,  U \text { a set }\}.$$
Here $SP$ denotes the operation of forming subdirect products. It is straightforward to show that ${\bf RPA}_{\alpha}\subseteq \PA_{\alpha}.$ 
Daigneault and Monk \cite{DM} proved that for $\alpha\geq \omega$ the converse inclusion also holds, 
that is ${\bf RPA}_{\alpha}=\PA_{\alpha}.$ 
This is a completeness theorem for certain infinitary extensions of first order logic 
without equality \cite{K}.

In this paper we are concerned with the following question:
If $\A$ is a polyadic algebra,  is there a 
representation of $\A$ that preserves infinite meets and joins, whenever they exist? (A representation of a given abstract algebra is basically a non-trival 
homomorphism from this algebra into a set algebra).
To make the problem more tangible we need a few preparations.
In what follows $\prod$ and $\sum$ denote infimum
and supremum, respectively. We will encounter situations where we need to evaluate a supremum of a 
given set in more than one algebra, in which case we will add a superscript to the supremum indicating the algebra we want. 
For set algebras, we identify notationally the algebra with its universe, since the operations 
are uniquely defined given the unit of the algebra.

Let $\A$ be a polyadic algebra and $f:\A\to \wp(^{\alpha}U)$ be a representation of $\A$.
If $s\in X$, we let
$$f^{-1}(s)=\{a\in \A: s\in f(a)\}.$$
An atomic representation $f:\A\to \wp(^{\alpha}U)$ is a representation such that for each 
$s\in V$, the ultrafilter $f^{-1}(s)$ is principal. 
A complete representation of $\A$ is a representation $f$ satisfying
$$f(\prod X)=\bigcap f[X]$$
whenever $X\subseteq \A$ and $\prod X$ is defined.

The following is proved by Hirsch and Hodkinson for cylindric algebras. The proof works for $\PA$'s without any modifications.

\begin{lemma}\label{r} Let $\A\in \PA_{\alpha}$. A representation $f$ of $\A$
is atomic if and only if it is complete. If $\A$ has a complete representation, then it is atomic.
\end{lemma}
\begin{demo}{Proof}\cite{HH}
\end{demo}
By Lemma \ref{r} a necessary condition for the existence of complete representations is the condition of atomicity.
We now prove a converse to this result, namely, that when $\A$ is atomic, then 
$\A$ is completely representable. 
We need to recall from \cite[definition~5.4.16]{HMT2},  the notion of neat reducts of polyadic algebras, 
which will play a key role in our proof of the main theorem. 
\begin{definition} Let $J\subseteq \beta$ and 
$\A=\langle A,+,\cdot ,-, 0, 1,{\sf c}_{(\Gamma)}, {\sf s}_{\tau}\rangle_{\Gamma\subseteq \beta ,\tau\in {}^{\beta}\beta}$
be a $\PA_{\beta}$.
Let $Nr_J\B=\{a\in A: {\sf c}_{(\beta\sim J)}a=a\}$. Then
$$\Nr_J\B=\langle Nr_{J}\B, +, \cdot, -, {\sf c}_{(\Gamma)}, {\sf s}'_{\tau}\rangle_{\Gamma\subseteq J, \tau\in {}^{\alpha}\alpha}$$
where ${\sf s}'_{\tau}={\sf s}_{\bar{\tau}}.$ Here $\bar{\tau}=\tau\cup Id_{\beta\sim \alpha}$.
The structure $\Nr_J\B$ is an algebra, called the {\it $J$ compression} of $\B$.
When $J=\alpha$, $\alpha$ an ordinal, then $\Nr_{\alpha}\B\in \PA_{\alpha}$ and is called the {\it neat $\alpha$ reduct} of $\B$ 
and its elements are called 
$\alpha$-dimensional.
\end{definition}

The notion of neat reducts is also extensively studied for cylindric algebras \cite{neat}.
We also need, \cite [theorem~ 2.1]{DM} and top of  p.161 in op.cit:

\begin{definition} Let $\A\in \PA_{\alpha}$. 
\begin{enumroman}

\item If $J\subseteq \alpha$, an element $a\in A$ is {\it independent of $J$} if ${\sf c}_{(J)}a=a$;
$J$ supports $a$ if $a$ is independent of $\alpha\sim J$.

\item The {\it effective degree} of $\A$ is the smallest cardinal $\mathfrak{e}$ such that each element of $\A$ admits a support 
whose cardinality does not exceed $\mathfrak{e}$.
\item The {\it local degree} of $\A$ is the smallest cardinal $\mathfrak{m}$ such that each element of $\A$ has support of cardinality $<\mathfrak{m}$. 
\item The {\it effective cardinality} of $\A$ is $\mathfrak{c}=|Nr_JA|$ where $|J|=\mathfrak{e}.$ (This is independent of $J$). 
\end{enumroman}
\end{definition} 

We shall use the following elementary known facts about boolean algebras and topological spaces. 
\begin{enumarab}
\item Let $\B$ be a boolean algebra, 
and let $S$ be its Stone space whose underlying set consists of all ultrafilters of $\B$. The topological space $S$ has a clopen base 
of sets of the form $N_b=\{F\in S: b\in F\}$ for $b\in B$. Assume that $X\subseteq B$ and $c\in B$
are such that $\sum X=c$. Then the set $N_c\sim \bigcup_{x\in X}N_x$ is nowhere dense in the 
Stone topology. In particular, if $c$ is the top element, then it follows that $S\sim \bigcup_{x\in X}N_x$ is nowhere dense.
(A nowhere dense set is one whose closure has empty interior).
\item Let ${\bf X}=(X, \tau)$ be a topological space. Let $x\in X$ be an isolated point in the sense that there is an open set $G\in \tau$ containing $x$, 
such that
$G\cap X=\{x\}$. Then $x$ cannot belong to any nowhere dense subset of $X$.  
\end{enumarab}
The proofs of these facts are entirely straightforward. They follow from the basic definitions.

Now we formulate and prove the main result. The proof is basically a Henkin construction together with a simple 
topological argument. The proof also has affinity with the proofs of the main theorems in \cite{DM} and \cite{super}.

\begin{theorem} Let $\alpha$ be an infinite ordinal. Let $\A\in \PA_{\alpha}$ be atomic. 
Then $\A$ has a complete representation.
\end{theorem}

\begin{demo}{Proof} Let $c\in A$ be non-zero. We will find a set $U$ and a homomorphism 
from $\A$ into the set algebra with universe $\wp(^{\alpha}U)$ that preserves arbitrary suprema
whenever they exist and also satisfies that  $f(c)\neq 0$. $U$ is called the base of the set algebra. 
Let $\mathfrak{m}$ be the local degree of $\A$, $\mathfrak{c}$ its effective cardinality and 
$\mathfrak{n}$ be any cardinal such that $\mathfrak{n}\geq \mathfrak{c}$ 
and $\sum_{s<\mathfrak{m}}\mathfrak{n}^s=\mathfrak{n}$. The cardinal 
$\mathfrak{n}$ will be the base of our desired representation.

Now there exists $\B\in \PA_{\mathfrak{n}}$ 
such that $\A\subseteq \Nr_{\alpha}\B$ and $A$ generates $\B$. 
The local degree of $\B$ is the same as that of $\A$, 
in particular each $x\in \B$ admits a support of cardinality $<\mathfrak{m}$. Furthermore, $|\mathfrak{n}\sim \alpha|=|\mathfrak{n}|$ and
for all $Y\subseteq A$, we have $\Sg^{\A}Y=\Nr_{\alpha}\Sg^{\B}Y.$ 
All this can be found in \cite{DM}, see the proof of theorem 1.6.1 therein; in such a proof, 
$\B$ is called a minimal dilation of $\A$. 
Without loss of generality, we assume that $\A=\Nr_{\alpha}\B$,  since $\Sg^{\A}A=\Nr_{\alpha}\Sg^{\B}A=\Nr_{\alpha}\B.$ (In the last equality 
we are using that $A$ generates $\B$).
Hence $\A$ is first order interpretable in $\B$. In particular, any first order sentence (e.g. the one expressing that $\A$ is atomic)
of the language of $\PA_{\alpha}$ translates effectively to a sentence $\hat{\sigma}$ of the language of $\PA_{\beta}$ such that for all 
$\C\in \PA_{\beta}$, we have
$\Nr_{\alpha}\C\models \sigma\longleftrightarrow \C\models \hat{\sigma}$. 
Here the languages are uncountable, even if $\A$ is countable and has countable dimension,  so we use effective in a loose sense, 
but it roughly means that there is an effective 
procedure or algorithm that does 
this translation. Since $\A=\Nr_{\alpha}\B$ and $\A$ is atomic, it 
follows that $\B$ is also atomic.
Let $\Gamma\subseteq \alpha$ and $p\in \A$. Then in $\B$ we have, see \cite{DM} the proof of theorem 1.6.1,  
\begin{equation}\label{tarek1}
\begin{split}
{\sf c}_{(\Gamma)}p=\sum\{{\sf s}_{\bar{\tau}}p: \tau\in {}^{\alpha}\mathfrak{n},\ \  \tau\upharpoonright \alpha\sim\Gamma=Id\}.
\end{split}
\end{equation}
Here, and elsewhere throughout the paper,  for a transformation $\tau$ wth domain $\alpha$ and range included in $\mathfrak{n}$,
$\bar{\tau}=\tau\cup Id_{\mathfrak{n}\sim \alpha}$. 
Let $X$ be the set of atoms of $\A$. Since $\A$ is atomic, then  $\sum^{\A} X=1$. By $\A=\Nr_{\alpha}\B$, we also have $\sum^{\B}X=1$.
We will further show that
for all $\tau\in {}^{\alpha}\mathfrak{n}$ we have,
\begin{equation}\label{tarek2}
\begin{split}
\sum {\sf s}_{\bar{\tau}}^{\B}X=1.
\end{split}
\end{equation}
It suffices to show that for every $\tau\in {}^{\alpha}\alpha$ and $a\neq 0\in \A$, there exists $x\in X$, such that ${\sf s}_{\tau}x\leq a.$
This will show that for any $\tau\in {}^{\alpha}\alpha$, the sum  $\sum{\sf s}_{\tau}^{\A}X$ 
is equal to the top element in $\A$, which is the same as that of $\B$. 
The required will then follow since for $Y\subseteq A$, we have
$\sum^{\A}Y=\sum ^{\B}Y$ by $\A=\Nr_{\alpha}\B$, and for all $x\in A$ and $\tau\in {}^{\alpha}\alpha$, we have 
${\sf s}_{\tau}^{\A}x={\sf s}_{\bar{\tau}}^{\B}x$, since $\B$ is a dilation of $\A$.

Our proof proceeds by certain non-trivial manipulations of substitutions. 
Assume that $\tau\in {}^{\alpha}\alpha$ and non-zero $a\in A$ are given. Suppose for the time being that $\tau$ is onto. 
We define a right inverse $\sigma$ of $\tau$ the usual way. That is for $i\in \alpha$ choose $j\in \tau^{-1}(i)$ and set $\sigma(i)=j$. 
Then $\sigma\in {}^{\alpha}\alpha$ and in fact $\sigma$ is one to one. 
Let $a'={\sf s}_{\sigma}a$. Then $a'\in A$ and  $a'\neq 0$, for if it did, then by polyadic axioms (7) and (8), we would get 
$0={\sf s}_{\tau}a'={\sf s}_{\tau\circ \sigma}a=a$ 
which is not the case, since $a\neq 0$.
Since $\A$ is atomic, then there exists an atom $x\in X$ such that $x\leq a'$.
Hence using that substitutions preserve the natural order on the boolean algebras in question, 
since they are boolean endomorphisms, and axioms (7) and (8)
in the polyadic axioms, we obtain
$${\sf s}_{\tau}x\leq {\sf s}_{\tau}a'={\sf s}_{\tau}{\sf s}_{\sigma}a={\sf s}_{\tau\circ \sigma}a={\sf s}_{Id}a=a.$$
Here all substitutions are evaluated in the algebra $\A$, and we are done in this case.
Now assume that $\tau$ is not onto. Here we use the spare dimensions of $\B$. 
By a simple cardinality argument, baring in mind that $|\mathfrak{n}\sim \alpha|=|\mathfrak{n}|$, 
we can find $\bar{\tau}\in {}^{\mathfrak{n}}\mathfrak{n}$, such that 
$\bar{\tau}\upharpoonright \alpha=\tau$, and 
$\bar{\tau}$ is onto. We can also assume that $\bar{\tau}\upharpoonright (\mathfrak{n}\sim \alpha)$ is one to one since 
$|\mathfrak{n}\sim \alpha|=|\mathfrak{n}\sim Range(\tau)|=\mathfrak{n}.$
Notice also that we have $\bar{\tau}^{-1}(\mathfrak{n}\sim \alpha)\cap \alpha=\emptyset$.  
Indeed if $x\in \bar{\tau}^{-1}(\mathfrak{n}\sim \alpha)\cap \alpha$, then $\bar{\tau}(x)\notin \alpha$, while $x\in \alpha$, 
which is impossible, since $\bar{\tau}\upharpoonright \alpha=\tau$, 
so that $\bar{\tau}(x)=\tau(x)\in \alpha$.
Now let $\sigma\in {}^{\mathfrak{n}}\mathfrak{n}$ be a right inverse of $\tau$ (on $\mathfrak{n}$), so that
for chosen $j\in \bar{\tau}^{-1}(i)$, we have $\sigma(i)=j.$
We distinguish between two cases:

(i) $\sigma(\alpha)\subseteq \alpha$.
Let $a'={\sf s}_{\sigma}^{\B}a={\sf s}_{\sigma\upharpoonright \alpha}^{\A}a.$ The last equation holds because $\B$ is a dilation of $\A$.
Then $a'\in A$, since $\sigma(\alpha)\subseteq \alpha$. Also $a'\neq 0$, by same reasoning as above.
Now there exists an atom $x\in X$, such that $x\leq a'$. A fairly straightforward computation, using that $\B$ is a dilation of $\A$,
together with the fact that substitutions preserve order, and axioms 
(7) and (8)
in the polyadic axioms applied to $\B$, gives
$${\sf s}_{{\tau}}^{\A}x\leq {\sf s}_{{\tau}}^{\A}a'={\sf s}_{\bar{\tau}}^{\B}a'= {\sf s}_{\bar{\tau}}^{\B}{\sf s}_{\sigma}^{\B}a=
{\sf s}_{\bar{\tau}\circ \sigma}^{\B}a={\sf s}_{Id}^{\B}a=a,$$
which finishes the proof in the first case.

(ii) $\sigma(\alpha)$ is not contained in $\alpha.$ 

We proceed as follows. Let $a'={\sf s}_{\sigma}^{\B}a\in \B$, then by the same reasoning as above, $a'\neq 0$. 
But $a'\in \B$, hence, there exists $\Gamma\subseteq \mathfrak{n}\sim \alpha$, such that 
$a''={\sf c}_{(\Gamma)}a'\in \A,$ since $\A=\Nr_{\alpha}\B$. But  $a'\leq a''$, so $a''$ is also a non-zero element in $\A$. 
Let $x\in X$ be an atom below $a''$. Then we have
$${\sf s}_{\tau}^{\A}x\leq {\sf s}_{\bar{\tau}}^{\B}a''={\sf s}_{\bar{\tau}}^{\B}{\sf c}_{(\Gamma)}^{\B}a'=
{\sf s}_{\bar{\tau}}^{\B}{\sf c}_{(\Gamma)}^{\B}{\sf s}_{\sigma}a=
{\sf c}_{(\Delta)}^{\B}{\sf s}_{\bar{\tau}}^{\B}{\sf s}_{\sigma}^{\B}a.$$
Here $\Delta=\bar{\tau}^{-1}\Gamma$, and the last equality follows from axiom (10) in the polyadic axioms noting that
$\bar{\tau}\upharpoonright \Delta$ is one to one since  
$\bar{\tau}\upharpoonright (\mathfrak{n}\sim \alpha)$ is one to one and  
$\Delta\subseteq \mathfrak{n}\sim \alpha$. The last inclusion holds by noting that $\bar{\tau}^{-1}(\mathfrak{n}\sim \alpha)\cap \alpha=\emptyset,$ 
$\Gamma\subseteq \mathfrak{n}\sim \alpha,$ and $\Delta=\bar{\tau}^{-1}\Gamma.$
But then going on, we have also by axioms (7) and (8) of the polyadic axioms
$${\sf c}_{(\Delta)}^{\B}{\sf s}_{\bar{\tau}}^{\B}{\sf s}_{\sigma}^{\B}a= 
{\sf c}_{(\Delta)}^{\B}{\sf s}_{\bar{\tau}\circ \sigma}^{\B}a={\sf c}_{(\Delta)}^{\B}{\sf s}_{Id}^{\B}a={\sf c}_{(\Delta)}^{\B}a=a.$$
The last equality follows from the fact that $\Delta\subseteq \mathfrak{n}\sim \alpha,$ 
and $a$ is $\alpha$-dimensional, that is $a\in \A=\Nr_{\alpha}\B.$  
We have proved that ${\sf s}_{\tau}^{\A}x\leq a$, and so we are done with the second slightly more difficult case, as well.

Let $S$ be the Stone space of $\B$, whose underlying set consists of all boolean ulltrafilters of 
$\B$. Let $X^*$ be the set of principal ultrafilters of $\B$ (those generated by the atoms).  
These are isolated points in the Stone topology, and they form a dense set in the Stone topology since $\B$ is atomic. 
So we have $X^*\cap T=\emptyset$ for every nowhere dense set $T$ (since principal ultrafilters, which are isolated points in the Stone topology,
lie outside nowhere dense sets). 
For $a\in \B$, let $N_a$ denote the set of all boolean ultrafilters containing $a$.
Now  for all $\Gamma\subseteq \alpha$, $p\in B$ and $\tau\in {}^{\alpha}\mathfrak{n}$, we have,
by the suprema, evaluated in (1) and (2):
\begin{equation}\label{tarek3}
\begin{split}
G_{\Gamma,p}=N_{{\sf c}_{(\Gamma)}p}\sim \bigcup_{{\tau}\in {}^{\alpha}\mathfrak{n}} N_{s_{\bar{\tau}}p}
\end{split}
\end{equation}
and
\begin{equation}\label{tarek4}
\begin{split}
G_{X, \tau}=S\sim \bigcup_{x\in X}N_{s_{\bar{\tau}}x}.
\end{split}
\end{equation}
are nowhere dense. 
Let $F$ be a principal ultrafilter of $S$ containing $c$. 
This is possible since $\B$ is atomic, so there is an atom $x$ below $c$; just take the 
ultrafilter generated by $x$.
Then $F\in X^*$, so $F\notin G_{\Gamma, p}$, $F\notin G_{X,\tau},$ 
for every $\Gamma\subseteq \alpha$, $p\in A$
and $\tau\in {}^{\alpha}\mathfrak{n}$.
Now define for $a\in A$
$$f(a)=\{\tau\in {}^{\alpha}\mathfrak{n}: {\sf s}_{\bar{\tau}}^{\B}a\in F\}.$$
Then $f$ is a homomorphism from $\A$ to the full
set algebra with unit $^{\alpha}\mathfrak{n}$, such that $f(c)\neq 0$. We have $f(c)\neq 0$ because $c\in F,$ so $Id\in f(c)$. 
The rest can be proved exactly as in \cite{super}; the preservation of the boolean operations and substitutions is fairly 
straightforward. Preservation of cylindrifications
is guaranteed by the condition that $F\notin G_{\Gamma,p}$ for all $\Gamma\subseteq \alpha$ and all $p\in A$. (Basically an elimination 
of cylindrifications, this condition is also used in \cite{DM}
to prove the main representation result for polyadic algebras.)
Moreover $f$ is an 
atomic representation since $F\notin G_{X,\tau}$ for every $\tau\in {}^{\alpha}\mathfrak{n}$, 
which means that for every $\tau\in {}^{\alpha}\mathfrak{n}$, 
there
exists $x\in X$, such that
${\sf s}_{\bar{\tau}}^{\B}x\in F$, and so $\bigcup_{x\in X}f(x)={}^{\alpha}\mathfrak{n}.$ 
We conclude that $f$ is a complete  representation by Lemma \ref{r}.
\end{demo}
Contrary to cylindric algebras, we have:
\begin{corollary} The class of completely representable polyadic algebras of infinite dimension is elementary
\end{corollary}
\begin{demo}{Proof} Atomicity can be expressed by a first order sentence.
\end{demo}
\section{ A metalogical reading; an extension of a theorem of Vaught}

Polyadic algebras of infinite dimension correspond 
to a certain infinitary logic studied by Keisler, and referred to in the literature as Keisler's logic. Keisler's 
logic allows formulas of infinite length and quantification on infinitely many 
variables, with semantics defined as expected.
While Keisler \cite{K}, and independently Monk and Daigneault \cite{DM}, 
proved a completeness theorem for such logics, our result implies a 
`Vaught's theorem' for such logics, namely, that atomic theories have atomic models, in a sense to be made precise.

Let $\L$ denote Keislers's logic with $\alpha$ many variables 
($\alpha$ an infinite ordinal). For a structure $\M$, a formula $\phi$,  and an assignment $s\in {}^{\alpha}M$, we write
$\M\models \phi[s]$ if $s$ satisfies $\phi$ in $\M$. We write $\phi^{\M}$  for the set of all assignments satisfying $\phi.$

\begin{definition} Let $T$ be a given $\L$ theory. 
\begin{enumarab}
\item A formula $\phi$ is said to be complete in $T$ iff for every formula $\psi$ exactly one of 
$$T\models \phi\to \psi, \\ T\models \phi\to \neg \psi$$
holds.
\item A formula $\theta$ is completable in $T$ iff there there is a complete formula $\phi$ with $T\models \phi\to \theta$.
\item $T$ is atomic iff if every formula consistent with $T$ is completable in $T.$

\item A model $\M$ of $T$ is atomic if for every $s\in {}^{\alpha}M$, there is a complete formula $\phi$ such that $\M\models \phi[s].$
\end{enumarab}
\end{definition}
We denote the set of formulas in a given language by $\Fm$ and for a set of formula $\Sigma$ we write $\Fm_{\Sigma}$ for the Tarski-
Lindenbaum quotient (polyadic)
algebra. \begin{theorem} Let $T$ be an atomic theory in $\L$ and assume that $\phi$ is consistent with $T$. Then $T$ has an atomic model in which
$\phi$ is satisfiable.
\end{theorem}
\begin{demo}{Proof}  Assume that $T$ and $\phi$ are given. Form the Lindenbaum Tarski algebra $\A=\Fm_T$ and let $a=\phi/T$.
Then $\A$ is an atomic polyadic algebra, since $T$ is atomic, and $a$ is non-zero, because $\phi$ is consistent with $T$.
Let $\B$ be a set algebra with base $M$, and $f:\A\to \B$ be a complete representation
such that $f(a)\neq 0$.
We extract a model $\M$ of $T$, with base $M$, from $\B$ as follows. 
For  a relation symbol $R$ and $s\in {}^{\alpha}M$, $s$ satisfies $R$ if $s\in h(R(x_0,x_1\ldots ..)/T)$. Here the variables occur in their 
natural order.
Then one can prove by a straightforward induction that $\phi^{\M}=h(\phi/T)$. Clearly $\phi$ is satisfiable in $\M$. 
Moreover, since the representation 
is complete it readily follows that  $\bigcup\{\phi^{\M}: \phi\text { is complete }\}={}^{\alpha}M$, and we are done.
\end{demo}

For ordinary first order logic atomic theories  in countable languages have atomic models, 
as indeed Vaught proved, but in the first order context countability is essentially 
needed. 

Our theorem can also be regarded as an omitting types theorem, for possibly uncountable languages, for the representation constructed in our theorem
omits the set (or infinitary type) $X^{-}=\{-x: x\in X\}$, in the sense that the representation $f$, defined in our main theorem, satisfies 
$\bigcap_{x\in X^{-}}f(x)=\emptyset$.  A standard omitting types theorem for Keisler's logic, addressing the omission of 
a family of types, 
not just one, is highly 
problematic since, even in the countable case, i.e when the universe of algebras is countable, since we have uncountably many 
operations. 

Nevertheless, a  natural omitting types theorem can be formulated as follows. 
Let $\L$ denote Keisler's logic, and let $T$ be an $\L$ theory. A set $\Gamma\subseteq \Fm$ is principal, if there exist a formula $\phi$ consistent with $T$,
such that $T\models \phi\to \psi$ for all $\psi\in \Gamma$. Otherwise $\Gamma$ is non-principal. A model $\M$ of $T$ omits $\Gamma$, if
$\bigcap_{\phi\in \Gamma}\phi^{\M}=\emptyset.$  Then the omitting types theorem in this context says:
If $\Gamma$ is non-principal, then there is a model $\M$ of $T$ that omits $\Gamma$.
Algebraically:

\begin{athm}{OTT} Let $\A\in \PA_{\alpha}$, and $a\in A$ be non-zero. Assume that $X\subseteq A$, is such that $\prod X=0$.
Then there exists a set algebra $\B$ 
and a homomorphism $f:\A\to \B$ such that $f(a)\neq 0$ and $\bigcap_{x\in X}f(x)=\emptyset.$
\end{athm}

Note that if $X$ satisfies $\prod {\sf s}_{\tau}X=0$, for evey $\tau\in {}^{\alpha}\alpha$,
then ${\bf OTT}$ holds (by the proof of our main theorem). We do not know whether ${\bf OTT}$ holds for arbitray $X$ with $\prod X=0$.

Unlike omitting types theorems for countable languages, our proof does not resort to the Baire Category theorem 
\cite{s}, \cite{t}. We believe that our result is an important addition to the model theory of Keisler's logics, a highly unexplored area. 
Beth definability and Craig interpolation for Keisler's logic are proved in \cite{super} 

\section{Concluding Remarks}

\begin{enumarab}

\item The technique used above is a central technique in the representation theory of both polyadic algebras \cite{super} and cylindric algebras 
\cite{neat}, \cite{s} together with  
a simple topological argument, namely, that isolated points in the Stone topology lie outside nowhere dense sets.
The former technique is called the neat embedding technique, which is an algebraic version of Henkin's completeness results, 
where the spare dimensions of $\B$ are no more than 
added witnesses that eliminate quantifiers in a sense.

\item The class of completely representable quasipolyadic {\it equality} algebras of infinite dimension, where only substitutions corresponding to 
finite transformations, is not elementary.
This can be proved exactly like the cylindric case \cite{HH} using a simple cardinality argument. However, 
the case of quasipolyadic algebras of infinite dimensions {\it without equality} remains an open 
problem. The finite dimensional case for such algebras is completely settled in \cite{Moh1}.

\item It is commonly accepted that polyadic algebras and cylindric algebras belong to different paradigms. For example, any universal axiomatization 
of the class of representable cylindric algebras of dimension $>2$  has an inevitable degree of complexity \cite{Andreka}, 
but  the class of representable polyadic algebras is  axiomatized by a fairly simple schema.
While the class of polyadic algebras has the superamalgamation property (which is a strong form of amalgamation) \cite{super}, 
the class of (representable)  cylindric algebras of dimension $>1$ fails to have even the amalgamation property \cite{amal}.
This paper manifests another dichotomy between those two paradigms, because the class of completely representable 
cylindric algebras  of dimension $>2$ is not elementary \cite{HH}, while the class of
completely representable polyadic algebras is. 

\item The case of polyadic agebras of infinite dimension {\it with equality} is much more involved. 
In this case the class of representable algebras is not a variety; it is not closed under ultraproducts, 
although every algebra has the neat embedding property (can be embedded into the neat reduct of algebras in every higher 
dimension).  In particular, we do not guarantee that atomic algebras are even representable, let alone admit a complete representation. 
Still we can ask whether atomic {\it representable} algebras are completey representable. 
The question seems to be a hard one, because we cannot resort to a neat embedding theorem as we did here for the equality free case.
\end{enumarab}

\end{document}